# HISTORIA MATHESEOS. EARLY DEVELOPMENT STAGE

## HISTORY OF MATHEMATICS. HISTORIOGRAPHY

Sinkevich G.I.

Saint Petersburg State University of Architecture and Civil Engineering Vtoraja Krasnoarmejskaja ul. 4, St. Petersburg, 190005, Russia

**Summary.** This article focuses on evolvement of the history of mathematics as a science and development of its methodology from the 4th century B.C. to the age of Enlightenment.

The history of mathematics was evolving together with the mathematics itself. Its methodology gradually established. Having started with short excursuses in the history of certain issues and biographical details of scientists, the history of mathematics ended up with research using historical, textual, and mathematical methods, and achieved significant results. One would find wealth of diverse literature devoted to the history of mathematics – from popular literature to rigorous research. People in all countries studied the history of mathematics; it is included in education courses; it is interesting for all lovers of mathematics.

## *4th century B.C., Eudemus of Rhodes.*

Mathematics was laid down as a science in Ancient Greece, and the first famous work devoted to the history of mathematics was the *History of Geometry* (Γεωμετρικὴ ἱστορία) by Eudemus of Rhodes[1], Aristotle's student. This work was repeated in *Comments to Euclid* by Proclus (The Catalogue of Geometers) [1, p. 82-86].

Eudemus is known to have also written *The History of Arithmetic* and *The History of Astronomy,* none of which has survived to our days and are only mentioned by other authors. He concluded that arithmetic resulted from trade activities of the Phoenicians, and geometry, from geodetic activities of Egyptians. He described a period of three centuries in Greek mathematics in his works. According to L.Y. Zhmud, researcher of the ancient world, Eudemus laid emphasis on three aspects: he specified the "ground breaker", looked for a faithful mathematical proof, and compared works of several scientists who studied the same issue [2, p. 277].

"The main question he raised in his Histories, was 'who discovered and what?'" [2, p. 289]. Eudemus arranged his materials chronologically – from discoveries of the teacher to discoveries of disciples.

Proclus began classifying mathematical literature by goal: "subject of investigation and for students" [1, p. 87]. Therefore, he compared teachers' merits as well."

Whereas the attitude to discovery of mathematical results was heuristic, meaning that ancient mathematicians 'discovered', not 'created', the history of mathematics was presented as a history of discoveries, i.e. as a linear flow 'who follows whom'.

Works of ancient mathematicians was normally preceded by a historical journey, for example, in Archimedes' works (3rd century B.C.). Starting with Diogenes Laertius (2nd–3rd century), biographical compendiums began to appear, i.e. summary versions of classic authors with comments which also included elements of history of mathematics. *Mathematical Collection* by Pappus of Alexandria (3rd–4th century) were a most complete compendium translated into Latin by Federico Commandino in 1588. It is thanks to Pappus that we now know many ancient problems. Let me also mention a compendium of Eutocius of Ascalon (5th–6th century) which contains solutions of the problem of cube duplication (Delosian problem) by ancient mathematicians and has preserved rendered fragments from *The History of Geometry* by Eudemus of Rhodes; and Simplicius' philosophical works ($5^{th}$–$6^{th}$ century). The genre of collections of rendered works (compendiums) became predominant in the Late Antiquity and Middle Ages up until the early Modern Age. Their subject was not only Greek mathematics, but Byzantine (e.g. treatises of a Byzantian, Maximus Planudes ($13^{th}$ century)), and a little of Indian mathematics as well.

---

[1] *Geometry and Geometers before Euclid* with a fragment from Eudemus of Rhodes was first published in Russian in a journal of
A.I. Goldenberg, "Mathematical Leaflet", v. I, 1879-1880.

Thabit ibn Qurra[2] (836–901), a mathematician and astronomer from Bagdade, translated works of Archimedes, Apollonius, Euclid, Ptolemy, and other ancient writers into Arabic. We can only read his translations of Archimedes' Treatises *Book of the Construction of the Circle Divided into Seven Equal Parts, On Tangent Circles*, as well as Books V–VII of *The Conics* by Apollonius. In Arabic mathematical literature of $8^{th}$–13th centuries, Arabic translations of ancient classic authors were accompanied by comments and supplemental information intended to systematise the knowledge of Greece, India, and the Arab world. Mathematicians' biographies were arranged chronologically in the Collections with extensive catalogues of manuscripts accompanied by historical comments. For example, in his *Introduction on Superiority of the Science of History* (*Prolegomena* or *Muqaddimah, Introduction to History*), an Arab historian Ibn Khaldun (1332–1406) reviewed the development of Greek and Islamic mathematics, dating its stages by periods of reign of tsars and khalifs. This book is available in French [3]. The treatise of Kâtip Çelebi, a Turkish historian of the $17^{th}$ century, *The Removal of Doubt from the Names of Books and the Arts,* containing 14,500 book titles and 10,000 names of authors and scholiasts in alphabetic order, which was published in Leipzig in seven volumes in Latin [4], is a most vivid example of bibliographic encyclopaedic genre.

As of the $12^{th}$ century, translations from Arab and Greek, or renderings of ancient mathematicians, and comments containing historical data, began to appear in Europe. Euclid was first translated or rendered by Adelard of Bath ($12^{th}$ century), an English scholastic philosopher who also translated al-Khwarizmi's astronomical tables and renderings of Euclid's Elements in 15 books by Giovanni Campano (Campanus, Campani, deceased in 1296). Campano himself did not know Arabic. He made his renderings with comments based on earlier translations.

With advent of book printing in late 15th century, Greek and Latin texts of Archimedes, Euclid, and other ancient authors started to appear. Scientists began to compare translations, analyse translators' and scholiasts' (relaters') errors.

## *1559, Johannes Buteo*

In 1559, a book of Johannes Buteo[3] (Buteo, Jean Borrel, 1492–between 1564 and 1572), *Quadrature of the Circle,* was published in two books where he invalidated many quadratures and sheltered Archimedes from hostile criticism; he also provided a list of errors made by Campani, Zamberti, Finé, and Peletier in their interpretations of Euclid in Latin [5]. Buteo analysed the errors of these translators of Euclid and Archimedes; provided approximate calculations of Bryson of Heraclea, Archimedes, and Ptolemy; and criticised the common misconception originating from Zamberti to the effect that the author of the demonstrations in Euclid's elements was Theon of Alexandria[4]. "Buteo has confidently mastered Archimedes' method and provided a summary of its use in the ancient world and Middle Ages." [6, p. 97].

---

[2] Sometimes, in Russian literature his name is spelled as Qurra.
[3] This is the Buteo which in 1559 calculated the capacity of Noah's Ark.
[4] Campano was the author of one of the first renderings of Euclid into Latin (*The Elements* in 15 books). Italian Zamberti (Zamberti, Zambertus, 1473 – after 1543) was the first to publish a printed translation of Euclid from Greek in 1505 (The Elements and other books). Zamberti corrected the mistakes in the medieval Campano's version in Latin. However, Zamberti was not a mathematician. Therefore, Luca Pacioli criticised him for his assaults upon Campano. In 1543, Tartaglia published his first translation of Euclid considering the text of Campano and Zamberti. Zamberti believed that Theon was the true author of the proof, while those were only the definitions and statements which belonged to Euclid. A French mathematician and cartographer Orontius (Oronce Finé, Orontius Finnæus, or Finæus, French Oronce Finé; 1494-1555) was Buteo's teacher. In 1532, he published his Protomathesis (Introduction in Mathematics) in Paris, where he explained the main notions used in The Elements by Euclid and the calculation of areas of planar figures as provided in Archimedes' works. He also described in this book his method of solving the circle quadrature problem which was later criticized by his student Johannes Buteo. Buteo criticized all the above authors including his teacher.

***1567, Pierre de la Ramée***. Buteo's book, as well as the book of Ramus (Pierre de la Ramée,1515–1572) contains no dates, adhering to the principle 'who follows whom'. In 1567, Ramée published his *Introduction to Mathematics divided into three books* [7] in Latin which formed a preamble to his large work, *Thirty-One Books of Mathematical Essays* [8]. This was a summary of prior discoveries in mathematics divided into three periods: from Adam to Abraham (Chaldean); from Abraham (Egyptian period); from Thales to Proclus and Theon of Alexandria (Greco-Roman period); and the fourth, modern period, from Theon ($5^{th}$ century) to Copernicus, Regiomontanus, and Cardano. The first book of the Introduction (p. 1–39) describes the first three periods and lists 65 names of Greek mathematicians. The second book – classification of mathematical sciences (Ramée described only arithmetic and geometry as a mathematical science; as to astronomy, optics, and music, he assigned them to physics) and their development in various European countries (p. 39–71). He listed around 30 names of $16^{th}$ century Reformers in his second book. They were mostly theologians, translators and scholiasts. Of mathematicians and astronomers, he mentioned, inter alia, Regiomontanus (however, Tycho Brahe, whom Ramée knew, was not mentioned); $16^{th}$ century: Herlinus, Copernicus, Rheticus (Copernicus' student), Rheinold, Santbecus, Leovitius (Cyprian Karasek Lvovicky), Dasypodius, Clavius, Landtgravius, Morshemius, Grunius, Xylander. Stiefel, however, was not mentioned. The third book described the development of teaching mathematical methods in European universities. The purpose of mathematics was described as practical application in trade, physics, architecture, astronomy, and other areas [8, p. 71–107]. The book contained almost no historical information – only rhetoric on the way ancient classic mathematicians should be presented in educational institutions. The book mentioned $16^{th}$ century mathematicians, such as Cardano, Maurolico, Piccolomini, Commandino, Tartaglia, and Dürer.

According to G.P. Matvievskaya, "In the fact that treatment of mathematics has not changed at all since Euclid, i.e. over two millennia, Ramée sees the proof of the "complexity and incomprehensibility" of its subject. On the other hand, this, of course, evidences of the clarity and perfection of The Elements. However, Ramus believed they must be submitted to a research regardless of the high standing of Euclid. As a result, he found a lot of material weaknesses in the renowned work. He believed, for example, that it was 'unmethodical' to begin giving an account of mathematics with geometry instead of algebra …, introduce definitions of mathematical notions before the need to define them arises" [9, p. 110].

G.N. Popov assessed the historical journey of Ramée as follows: "Given the deficiency of sources in the two first periods, the author need not speak much about them. However, the Greek science together with the Roman one was presented but on 36 pages. Dating is missing, although Ramus adheres to the historical sequence; he does not mention many geometers, which can of course be explained by the lack of sources he had at his disposal. However, those he used, evidence of the intelligent choice judging from the reliability of the obtained information." [10, p. 146]. Ramée noted that the copies of ancient manuscripts were preserved in Florence thanks to the Medici family[5], which evidences of his good knowledge of Italian historical literature. He set forth his speculations on the changes in methods of teaching mathematics in Christian Europe, giving preference to teaching in Protestant universities and criticising Aristotle. He mentioned Latin translations of Euclid and dissemination of information on Greek mathematics in Christian Europe. However, he mentioned nothing of the development of mathematics in the Islamic East. As in early works devoted to history of mathematics, Ramée considered mathematics as a combination of ancient Greek achievements remaining unchanged until the $17^{th}$ century to be looked up to, sometimes criticising teaching methods. He told nothing about Kepler's, Cardano's, or Tartaglia's results, although mentioned their names. There was no evolution dynamics of mathematics or its contents in Ramée's. In the emerging Age of Enlightenment, the notion of the progress of science had not yet become a historical category.

---

[5] Ramus' book was devoted to his patroness Ekaterina Medici, the Queen of France. However, it did not help to retrieve the Huguenot from ruin on the St. Bartholomew's night.

## New chronology

Chronology was an important problem in the history of mathematics as well as in general history. Each culture had its own chronology, and history of each culture did not correlate with other cultures; social time was described in each culture regardless of others. In Greece, they dated events by Olympiad; in the Arab world, from Hegira and by khalif; in Rome, time keeping was 'from the founding of Rome' (ab urbe condita); in Byzantium, 'from Adam', i.e. 'since the creation of the world'. Social time in different cultures was autonomous. The errors were significant. In 1582, the Pope, Gregory XIII published a bulla, *Inter gravissimas,* with an invitation to switch to a new calendar. First, some Catholic countries began gradually switching to the Gregorian calendar; thereafter, over a period of $17^{th}$ – $18^{th}$ centuries, Protestant countries, including Great Britain in 1752, were switching to it. Russia switched to the new calendar in 1918. It was Dionysius Exiguus who in the $6^{th}$ century suggested keeping time 'from the year of our Lord' (ab Anno Christi, ab inscriptione, Anno Domini). In Europe, this way of time keeping spread in late Middle Ages. In 725, in addition to the time keeping by Olympiad and by emperor, Bede the Venerable introduced absolute chronology for the first time. Beginning Chapter Two of Book One of his *Ecclesiastic History of the English People* (Historia ecclesiastica gentis anglorum), Bede wrote in 731: "ante incarnationis dominicae tempus" (before the incarnation of our Lord). This was the first time ever that the time countdown was mentioned. This is not to say that the countdown scale – B.C. and A.D. – developed in a prompt and natural way. Before the $16^{th}$ century, along with 'A.D.' time keeping, they used 'anno mundi' and many other ways of reckoning.

In the period from 1583 to 1629, books of Joseph Justus Scaliger (1540–1609) were published. He was a connoisseur of ancient culture and ancient calendars and the founder of modern chronology as an auxiliary science of history. Scaliger found ways of conversion between the systems of Ancient Rome, Ancient Greece, East Asia, and Mexico, using the method of astronomic dating of events by eclipse. This enabled him to correlate scientific discoveries in various cultures in the course of time.

In 1627, Dionisius Petavius (Denis Pétau, 1583–1652), a French scholar, suggested a 'before Christ' (ante Christum, B.C., century before Christian Era) system of counting down dates[6]. This system was universally recognized by the end of the $18^{th}$ century.

## Bernardino Baldi.

Bernardino Baldi (1553–1617), an Italian poet and mathematician, Commandino's student, was the first to try and use the new chronology in combination with former traditions of presenting history as a chronicle. It took him 12 years to create a *Chronicle of Mathematicians including their curriculum vitae* [11] as basis for a more comprehensive work. Among other preparatory materials for this work, his writings about Pythagor, Ctesibius, Hero of Alexandria, and Copernicus were preserved. The *Chronicle of Mathematicians* was written in Italian and contained around 200 life histories and an index of names. The book was written as a popular Who's Who. Names were arranged in a chronological order; in the left field, Greek chronology by Olympiad; in the right field, years before Christ or after Christ. E.g. about Euclid: 122 ($122^{nd}$ Olympiad) in the left field, 290 (anni avanti Christo – 290 before Christ) in the right field; and a text as follows: *"Euclid*. There is some evidence that the most esteemed mathematician from the City of Gela in Sicily studied in Alexandria and probably in Athens. We has written a lot, that is to say a book entitled *The Elements of Geometry* in which he surpassed all those who wrote before him, and he was so glorious that he was named στοιχειωτής – Stichiota[7](Sic! –G.S.). In addition to the *Elements*, he wrote a book entitled *Data*, *Porisms* in three volumes about perspective projection, about mirrors, a book about phenomena, and, to the best of my knowledge, a book entitled *Conics*, instead of the book about basics of music erroneously assigned to him. There is another apocryphal book about division of surfaces assigned to him by Mohammed from Bagdad. There is also a Plato section of Euclid which, according to Proclus, prepares the use of *The Elements* for the purposes of plotting Platonic solids." [11, p. 22-23].

---

[6] It should be noted that it's in these years that a new understanding of a number scale comes into existence too: minus (negative) numbers, as numbers that are less than zero (according to Stifel), are located to the left of (behind) zero. In 1629, A. Girard wrote about negative solutions of equations: "Solution using minus is explained in geometry as reversion, and minus retreats where plus goes further." [25, p. 228].
[7] The Author of "The Elements", Euclid the stoicheiotes.

Baldi also mentioned Arab and Persian mathematicians, as well as mathematicians from North Europe; many astronomers; some philosophers and theologians (e.g. Thomas Bradwardine, Nicholas of Cusa, Abraham Zacuto). Notably, there was no article devoted to Girolamo Cardano, although he mentioned his name in the articles devoted to Swineshead and Tartaglia. The timeframe covered by this book was from 545 B.C. to 1596 A.D.

## *1615, Joseph Blancanus*

In 1615, in his *Thesis on the Nature of Mathematics,* Joseph Blancanus (Giuseppe Biancani, 1566–1624), Italian mathematician and astronomer, represented the history of European and Asian mathematics as a history of discoveries in accordance with the new chronology, *Chronology of famous mathematicians* [12], in Latin. Although it contained certain inaccuracies[8], it was more complete than Ramée's work and included Islamic mathematicians. This was an attempt to use a unified time scale representing the European and Arab history of mathematical discoveries.

***1650, Gerardus Johann Voss.*** The book which enriched the history of mathematics with new methods was not written by a mathematician. It was written by a Dutch historian and philologist, Gerardus Johann Voss (Vossius, 1577–1649). His materials on the history of literature were such extensive that they included information on the development of mathematics as well. He collected these materials in his work *On the Nature and Structure of All Mathematical Sciences Supplemented by Mathematicians' Chronology* published posthumously [13] and republished as a part of the book entitled *On Four Main Arts, On Philology and Mathematical Sciences, Supplemented by Mathematicians' Chronology, issued in three books* [14]. Voss was not a mathematician and, at times, used inaccurate information, for which he was reasonably blamed by researchers. However, he was the first to use philological and source study methods in his historical and mathematical review.

Voss began with the history of alphabetical and digital numbers and symbols, systematised his presentation by section (geometry, arithmetic, optics, music, mechanics, logistics, geodesy, astronomy, calendar, chronology). Having mentioned Greek mathematicians, among others, he mentioned such mathematicians as Boethius, Alcuin, al-Farghani, Ibn al-Haytham, Sacrobosco, Nicholas of Cusa, Regiomontanus, Zakuto, Dürer, Copernicus, Maurolico, Cardano, Gemma, Commandino, Mercator, Ramée, Clavius, Viète, Ludolph van Ceulen, Tycho Brahe, Neper, van Roomen, Grégoire de Saint-Vincent, Stifel, Mersenne, Snellius, J. Golius, Cavalieri. He placed emphasis on translations of Greek classic authors into Arabic and thereafter, from Arabic into Latin [14, p. 55]; addressed works of Arabian historians. In the section devoted to the history of Alfonsine tables, Voss mentioned the notion of a progress of science for the first time – "progress of astronomy after Greeks" [14, p. 146]. For him, publishing a book, including a translation, or annotation, was a mathematical event. Voss' book was followed by an index rerum & verborum, i.e. an index of objects and words, an index of names and subjects with page numbers, and in addition, a list of printing mistakes. All this suggested a new sample form of a historical and mathematical research.

## *1674, C.-F. de Chales*

Claude-Francois de Chales, French mathematician, Jesuit professor (de Chales, Dechales, 1621–1678)[9] was the first in historiography to express his consciousness of the advance of mathematics in his *Treatise on the Advance of Mathematics and on Famous Mathematicians* which made part of Volume One of his three-volume edition of *The Course and World of Mathematics* [15], an pansophy which contained information from mathematics, physics, astronomy, astrology, and architecture. De Chales translated Euclid, and this translation was popular in France although it was worse than Roberval's translation. D. Smith wrote that, although de Chales published Euclid, his own contribution in the subject was more than modest. [16, p. 386].

---
[8] E.g. Thābit ibn Qurra (836-901) was described as a 13th century scientist; Roger Bacon (13th century), as a 14th century scientist; Leonardus Pisanus (Fibonacci, early 13th century), as a 15th century scientist

[9] Not to be confused with geometer Michel Chasles (1793 – 1880), French mathematician (geometer) and historian of mathematics

### *1681, J. Mabillon*

In 1681, a book of J. Mabillon (Jean Mabillon, 1632–1707) was published. Historian, originator of paleography, historical criticism, and chronology, wrote a book entitled *Diplomacy in Six Books* [17]. Mabillon understood diplomacy as a science addressing historical documents, evidence of their accuracy, methods of identifying forgery, ancient written instruments and materials, styles. Mabillon's book contained engraving plates with examples of ancient writing. Voss' and Mabillon's works influenced subsequent researchers and J. Wallis, in the first place.

### *1685, John Wallis*

His *Treatise of algebra both Historical and Practical. Shewing, the Original, Progress and Advancement thereof, from time to time; and by what Steps it hath attained to the Heighth at which now it is. With some additional Treatises. London. M.DC.LXXXV* [18] was published in 1685. Historians of mathematics (Cajori, Bobynin, Popov) blamed him of nationalism and lack of personal modesty which consisted in attributing discoveries of other mathematicians to himself (or his compatriots). This was a just reproach which, however, did not make his treatise less interesting, as it was written by an outstanding mathematician. The history of mathematics (algebra) here was for the first time presented as a history of ideas.

Speaking of works of ancient classic authors, Wallis gave names of translators and publishers. But it appears from the following example that the Arab history (chronology) had still existed apart from the European history. Wallis believed that Arabs knew algebra: "After Diophantus (if not before, also) this Learning was pursued by *Arabic* Authors (but little known in *Europe* for a long time)… Divers writers (is said) there are of *Algebra* in that Language, and from them (I suppose) the Denominations of *Diophantus* (if from him they learned it) came to be changed; and (beside the Denominations of Root, Square, and Cube) that of Sursolide (first, second, thirds, etc.) introduced. But I rather think the *Arabs*, either of themselves, or from some others, had it long before *Diophantus*, and think this reckoning of *Powers* (by Sursolids, etc.) different from Diophantus" [18, p. 5, Italics of the original]. According to Wallis, in England, algebra started to develop in the $12^{th}$ and $13^{th}$ centuries – earlier than in Europe – thanks to the fact that English scholastics knew Arabic. Englishmen used to visit Spain and bring many mathematical manuscripts with them. For example, in 1180, J. Morley (Morlacus, Morley, around 1140 – around 1210), mathematician and astronomer, studied Arab mathematical manuscripts in Toledo and brought a valuable set to England. The English were the first to translate Greek mathematical texts from Arabic. For example, in 1130, Adelard was the first to translate Euclid's Elements. Wallis mentioned an English theologian and historian, monk Bede Venerable (Saint Beda, Beda Venerabilis, late $7^{th}$ – early $13^{th}$ century), who wrote the history of English people; and then Alcuin (Alcuinus, around 735–804). Wallis erroneously called him Bede's student[10]. Wallis told in detail which Arab translations of ancient authors were brought to Oxford (including Merton College) and translated into Latin and English. In fact, this was a rendering of Voss' history of mathematics in the context of English history. He further told about the numerical values which originated from Moors and Arabs (p. 7) and Maximus Planudes. He also considered other number notations – Roman and Greek literal and number notations worldwide. He admitted that, although Arabic figures came from Saracens and Arabs, they originated from India. He compared Sacrobosco's (Johannes Sacrobosco, Sacrobosco, John of Holywood, around 1195–around 1256) numbering, who described fundamentals of Indo-Arabian numbering and arithmetic in his treatise *Algorithm* (*Algorismus de integris*): operations of addition, subtraction, averaging, duplication, multiplication, division, summing up arithmetical progressions, rooting, and cube-root extraction. Wallis believed that it was Luca Pacioli who was the first to bring the new notation to Europe (Luca de Burgo, Luca Pacioli, Summa de Arithmetica, Geometria, Proportioni, et Proportionalita, 1494). Thereafter, Wallis divided his narration into topics, which tells the difference between him and his predecessors.

---
[10] Alcuin was born after Bede had died, and was in tutelage of Archbishop Egbert, Bede's student

Further Wallis' synopsis.

Astronomic tables: Ptolemy, Copernicus. Decimal fractions appeared; then, logarithms. Archimedes' methods, including the method of using big numbers (with the help of 60-ary fractions, e.g. $1/4 = 15'$). Operations on them. On decimal fractions and use thereof in certain branches of arithmetic; on antiquity of decimal fractions. Described works of Briggs, Oughtred[11], Gellibrand[12] (*Trigonometria Britannica*, 1633), Regiomontan (1464), Ramus (1560), Schoener[13] (1585), Record (1550), Stevin (1585). Reducing fractions, or proportions, to a smaller number of characters with the nearest approximation to the real value: the development line from Archimedes: Van Gulen, Snellius[14]. Chapter 11: application of the same in relation to the proportion of the diameter and length of circumference (from Archimedes).

Chapter 12: *On Logarithms*. Wallis wrote about Neper, Briggs, Kepler, Rudolphian tables (1627), Mercator's *Logarithmotechnia* (1668), however, he did not mention logarithmic tables of J. Speidell[15] [19, 20, 21, p. 43] or the logarithmic rule invented by the English (astronomers Gunter and Wingate, mathematician Oughtred).

In Chapter 13, Wallis spoke about algebra: {…} [18, p. 61]. About Leonardo Pisano, Luca Pacioli, Cardano, Tartaglia, Nunes, Bombelli, and other authors of *Algebra* before Viète. He told a lot of good things about Pacioli (Pacioli's own books and his translation of Euclid)[16]. In Part Five of his *Sum*, Pacioli provided the basic materials on arithmetic as provided by ancient authors and his contemporaries. Subsequently, Wallis rendered Voss' story of Leonardo Pisano. Lucas de Burgo was the first to describe the abacus method, i.e. that of Luca Pacioli. In page 62, he mentioned Stifel, his *Arithmetic* of 1544, Rudolph, Cardano's *Arithmetic* and his *Great Art* (*Ars magna*) of 1545. Cardano's rule of solving a cubic equation which, according to Cardano, was found by Tartaglia too. Cardano's student Luigi Ferrari also added a bit[17]. In 1567, Pedro Nunes published *Algebra* in Spanish[18]. In 1579, Rafael Bombelli published the treatise *Algebra*[19] in Italian. In this book, he published the rule of solving a cubic and biquadratic equation as Tartaglia and Cardano had done before [18, p. 63].

We believe that Bobynin's observation to the effect that Wallis took the credit for the solution of the irreducible case of a cubic equation [10, p. 149] is groundless. On pages 173-174, Wallis spoke of the similarity of Cardano's and Herriot's methods. Further, in Chapters 46-48 (p. 175-181), he presented his own method of cube-root extraction from a binomial in the form of $a + \sqrt{b}$, where $b$ may have any sign. He found this method in 1647-48, and it was really similar to Bombelli's method published in 1572, and Wallis was familiar with Bombelli's book.

Then, Pierre de la Ramée (1570) was mentioned. The book was published by Schoener who also wrote books (*Numerical Geometry*). Leonard Digges[20] was one of our (i.e. Wallis') compatriots who wrote a book in 1579 entitled *Stratioticos* (military message deliverer). Another was Robert Recorde, 1552[21]. Chapter 14 was devoted to François Viète and his *symbolic arithmetic*[22].

Chapters 15 to 29 were devoted to Oughtred (Wallis was his disciple) and his book *The Key to Mathematics* (*Clavis Mathematicae*) [22], a textbook of arithmetic which was republished three times in Oughtred's lifetime and was thereafter used even in the 18th century. Wallis rendered it in the tiniest detail. Beginning with Chapter 30, Wallis wrote about Herriot's *Algebra*

---

[11] William (Guilelm) Oughtred, 1575–1660
[12] Henry Gellibrand, 1597–1636, Professor of Astronomy in Oxford who completed Briggs' unfinished work
[13] Lasarus Schonerus, Schoener, Schöner, 1543 – 1607, Ramée's publisher, scholiast, and partly a contributing author. He taught mathematics in Neustadt; was Provoste in Schmalkalden, Thuringia; and taught mathematics in Corbach high school
[14] Willebrord Snellius,1580–1626, Dutch mathematician and astronomer
[15] John Speidell, 1600–1634, teacher of mathematics in London, drafter of logarithmic tables
[16] Wallis called this period the Italian period. However, in effect, this was a revised publication of the Latin translation of Euclid made by Campano, where Pacioli had corrected numerous errors
[17] This is what Wallis wrote. However, this 'a bit' was a solution of a 4th-degree equation
[18] Wallis was mistaken here. It was in Portuguese
[19] This is the second publication. The first one was published in 1572
[20] Leonard Digges (around 1515 – around 1559), an English mathematician and topographer
[21] "If I be not misinformed" - Wallis' note
[22] Wallis meant *Francisci Vietae-in artem analyticem isagoge*.

[23], rendered it in detail, and asserted that Descartes adhered to Herriot. In this book, Herriot showed the way algebraic equations were laid down by multiplying linear binomials for the first time. Herriot rule (as Euler set it forth with reference to Herriot[23]) was as follows: each equation has as many positive roots as many variations of sign it contains, or as many negative roots as many repeat signs it contains. This only applies to those equations in which all roots are real [24, p. 468].

Now we refer to this rule as Descartes rule. Descartes himself claimed that, although he had Herriot's book (1631) at home, he only read it after he had finished his *Geometry* (1637). Descartes wrote to Constantijn Huygens (his father) in 1638: "Sir, I address my own books so seldom that one of your books turned to have been lost among them– although there is a mere half a dozen of them – and remained unnoticed for more than six months. It was Henriotti… I wanted to see this book, as they told me that it contained some calculation for geometry which was quite similar to mine; I found that this was correct, however, he did not go to the heart of the matter much and taught so few things on a plenty of pages, that there was no reason for me to have concerns about his thoughts which preceded mine" [25, p. 227. Also, see 26, p. 36-38]. In Chapter 53, Wallis accused Descartes of borrowings from Herriot. In particular, Wallis highly appreciated the innovation of Herriot who suggested that all members of an equation be written as one member of the equality, setting them to zero [24].

Wallis explained Herriot's seniority as follows:

"It will not be amiss here to insert a short Story which Dr. *John Pell*[25] lately told me he had from Sir *Charles Cavendish*[26], only Brother to *William*, then *Earl*, since First *Duce* of *Newcastle*; a Person of Honour, (well skilled in the Mathematics), who about that time lived in *Paris*. He discoursing there with Monsignor *Roberval*, concerning that piece of *Des Cartes* then lately published: I admire (saith *M. Roberval*) that notion in *Des Cartes* of putting over the whole Equation to one side, making it equal to Nothing, and how he lighted upon it. The reason why you admire it (saith Sir *Charles*) is because you are a French-man; for if you were an English-man, you would not admire it. Because (saith Sir *Charles*) we in England know whence he had it; namely from *Harriot's Algebra*. What book is that? (saith M. *Roberval*), I never saw it. Next time you come to my Chamber (saith Sir *Charles*) I will shew it you. Which a while after, he did: And upon perusal of it M. *Roberval* exclaimed with Admiration (Il l' aveu! Il l'aveu!) He had seen it! He had seen it! Finding all that in *Harriot* which he had before admired in *Des Cartes*; and not doubting that *Des Cartes* had it from thence. The Improvements of *Algebra* to be found in *Harriot* (as appears from what is already said), and which (all or most of them) we owe to him; (or which it will not be amiss, before I leave him, to give a brief Recapitulation) are chiefly these" [18, p. 198, Italics of the original].

Wallis was right about the seniority of Herriot. Descartes had Herriot's treatise when he was writing his *Geometry*. Many ideas which were thoroughly set forth and systematised by Descartes, had been first uttered by Herriot. It should be also noted that a book of A. Girard, *Invention nouvelle en Algèbre*, was published in Antwerp in 1629 to formulate the basic theorem in algebra eight years before Descartes. However, Wallis did not mention this fact. Further, in Chapter 55, [18, p. 208], Wallis repeated that Descartes' reasoning was based on Herriot's *Algebra* published in 1631, while Descartes' *Geometry* was published in 1637 in French and thereafter, in 1649 and in 1659, in Latin. Wallis demonstrated that, although used by Viète and Bombelli, many procedures could be asserted much easier based on Herriot's *Algebra*. This, certainly, does not derogate Descartes' role, who, unlike the English, did not proceed from geometry to algebra. Instead, developing algebra and generalizing the notion of a number, he set the analytical direction in the development of geometry.

---

[23] Wallis' book was in Euler's library in St. Petersburg

[24] It should be noted that Herriot did not select this form of notation as the final one. In his records, the constant term is more often on the right side. One can come across an equation written with a zero in the right side as far back as in Stifel's works.

[25] John Pell (1611–1685) was an English mathematician-algebraist..

[26] Sir Charles Cavendish (ca. 1594 – 1654) was an English aristocrat, Member of Parliament, and patron of philosophers and mathematicians. Cavendish knew Pell from the Welbeck period, along with the mathematicians Walter Warner and Robert Payne. He supported William Oughtred and knew John Wallis. Because the Cavendishes were royalist émigrés of the 1640s, the centre of this circle moved to Paris, where it took on the form of a salon. It grew around Thomas Hobbes and John Pell.

Wallis listed Herriot's achievements (up to 200-th page): symbols, terms, generation of equations by multiplying binomials, rule of signs (the number of positive and negative roots), methods of determining the number of real and imaginary roots, research of a quadratic equation, dividing an equation by a binomial, simplifying a cubic equation. Wallis acknowledged that almost all of these discoveries were made by Herriot, although some of them had been discovered earlier by Viète.

Wallis highly appreciated the role of Leonardo Pisano who reproduced Arab rules and symbols without resort to Diophantus who remained unknown in Europe until 1572[27]. According to Wallis, Stiefel was a good author who had never moved beyond quadratic equations [18, sheet a3]. Scipione del Ferro, Cardano, Tartaglia, and other developed a solution of a cubic equation. Bombelli took it a step further, solving biquadratic equations (with the help of cubic ones[28], reducing them to two quadratic equations). Nunes, Ramus, Schoener, Salignac[29], Clavius, Record, T. Digges[30], and some of our men (i.e. Englishmen – G.S.), were developing this subject in the last century. However, by and large, they had failed to take it a step further than the quadratic equations. At the same time, thanks to Xilander[31], Diophantus became known in Latin, and thereafter, thanks to Bachet, in Greek and Latin[32]; all his methods differed from Arab methods (followed by others). In particular, the procedure for naming exponents [18, sheet a3 verso]: using new symbols and figures were an important step in algebra. The next major step in the development of algebra was made by François Viète in 1590 in his *Specious Arithmetick*[33]. Wallis expressed his appreciation of this step. He noted that, unlike the preceding authors, in designating exponents, Viète adhered to Diophantus, not to Arabs.

Wallis' historical sketch was structured by problem. As a prominent mathematician, Wallis was by far more knowledgeable about mathematical information than other authors. He gave an unbiased account of the development of methods and discoveries in algebra and emerging analysis, although sometimes, one could, of course, blame him in subjective assessment. Unlike preceding authors, his narration was not a discrete set of biographies or discoveries. He showed mathematics as a continuous development of ideas and algebra in the first place. He demonstrated its internal relations and their continuity, genesis of mathematical knowledge, and creativity of a mathematician, not heuristicity. He discerned algebraic and geometrical methods and distinguished the inception of an analytical method, i.e. *Differential Calculation* method, in works of his contemporaries and, first of all, of Newton.

*1704, E. Bernard*. Edward Bernard (1638–1697) was a Savilian Professor[34] of astronomy in Oxford. He was connoisseur of ancient manuscripts; studied a lot of manuscripts of Apollonius of Perga; worked in Bodleian Library (Oxford) with Arab manuscripts brought from Spain, Morocco, Syria, Arab countries, and Turkey, which, to a large extent, replenished its collection. Edward Bernard found an Arab text of Apollonius entitled *Determined Section*; tried to recover those fragments that had been lost and translate them into Latin; he edited Josephus Flavius. Most of Edward Bernard's work consisted of annotating books from Bodleian Library: his work *On Ancient Weights and Measures (De mensuris et ponderibus antiquis* (1688)) was enclosed with a work of E. Pococke (1604–1691), an Orientalist scholar from Oxford. *Bernard's Catalogue* [27] comprised manuscripts from British and Irish libraries and served as a basic tool of scientists of that time. Many works of Bernard were not completed, which made his colleagues joke[35]. After Bernard 's death, his colleagues published a book about him [28, Section 9, p. 1-78] which included Bernard's work entitled *A Short List of Ancient Greek, Latin, and Arab Mathematicians, prepared by Dr. Eduardo Bernardo, the most honoured and educated man* [29, Section 11, p. 1-44] – annotated plan of republishing oeuvres of classic authors kept in European archives and libraries, on circa 44 pages. The adaptations and translations of Apollonius' *Conics*[36] Bernard had made were subsequently used by Edmond Halley (1656–1742) in the 1710 publication of Apollonius' works.

---

[27] Bombelli found a manuscript of Diophantus in the library of Vatican and published 143 problems in his *Algebra*. Wallis was mistaken about symbols. Leonardus Pisanus had no symbols. Wallis had not seen Pisanus' works. He learnt about them in Pacioli's works. (Thanks to J. Cesiano for this remark)

[28] This is not true! It was Ferrari who created a formula to solve a biquadratic equation! However, Bombelli uttered nothing about Ferrari, although Ferreri's formula had been set forth in Cardanus' *Ars magna* – G.S.

[29] Johannes Salignacus, Scottish

[30] Thomas Digges (1546 –1595), son of Leonardo Digges, English mathematician and astronomer, one of the first partisans and promoters of the heliocentric world

[31] Xilander published *Diophanti Alexandrini Rerum Arithmeticarum libri sex* in 1575 in Basel.

[32] Bachet de Meseriac published *Diophanti Alexandrini Arithmeticorum libri sex; et de Numeris multangulis liber unus. Nunc primum graece et latine editi, atque absolutissimis commentariis illustrati* in1621 in Paris.

[33] So Wallis calls Viiet's "species logistic" this way.

[34] In 1619, Sir Henry Savile, mathematician, custodian of Merton College in Oxford and provost of Eton College, complaining of the "poor condition of mathematical research in England", constituted two positions to be funded at his own expense: professor of geometry and professor of astronomy, which are in existence to the present day. The first professor of geometry was Henry Briggs

***1715, J. Raphson***. In 1715, a small posthumous edition of Newton's disciple, Joseph Raphson[37], was published. It was *History of Fluxions* [30], and the goal of this publication was to assert Newton's seniority in the discovery of differential calculation. Newton allowed Raphson to look through his works and his correspondence with Leibnitz, and the respective presentment of this correspondence in Raphson's book provided a strong support to Newton's position in this dispute.

***1741. Ch. Wolff***. Christian von Wolff (1679–1754), German philosopher, lawyer, and professor of mathematics, published a *Report on additions to mathematical sciences over one century* in Halley in 1707; *Mathematical Vocabulary* in Leipzig in 1716 [31] – a dictionary of mathematics in German, the best of those available by that time although not the first one, on 788 pages, the list of sources alone was on 54 pages; and an article entitled *A Summary of the Most Renowned Mathematical Works* in Volume V of *Elementary Fundamentals of Mathematical Sciences* [32, p. 3-168]. In Chapter One (p. 5-28), Wolff reviewed books, beginning with Euclid and finishing with publications of Academia Petropolitana until 1731. He, inter alia, mentioned works of young Euler. He devoted a paragraph with a summary to each of the above books. Those were books of French, English, and Dutch authors. Chapter Two (p. 29-32) was devoted to the history of arithmetic from Nicomachus to Neper. Chapter Three (p. 32-50), Geometry: Euclid and his translators, publishers, and annotators, European geometers, finishing with the year 1699. Chapter Four (51-69): analytical works from the ancient world to the inception of differential calculation. Much prominence was given to the dispute regarding the seniority of Newton and Leibnitz. And Wolff, professional lawyer and Leibnitz' friend, gave preference to Leibnitz in this matter, arranging their correspondence and publications chronologically, without getting involved in the mathematical substance. Chapter Five (p. 71-77): trigonometry from Ptolemy to Ozanam. Chapters Six to Thirteen were devoted to statics, mechanics (up to Euler), hydrostatics, aerometry, hydraulics, optics, catoptric, dioptric, perspective projection, astronomy, chronology, geography, gnomonics, civil architecture, pyrotechnics, and military architecture – well-established sections of the 18$^{th}$ century mathematics. The book contains an index of names.

***1742, I.H. Heilbronner***. The last book in the early period of historiography (before Montucl) was published in 1742. This was *The History of Mathematics at Large – from the Creation to the 16$^{th}$ Century A.D. including life stories of famous mathematicians, their doctrines, works, and manuscripts; in addition, a summary of main mathematical collections and works, and history of arithmetic to the present day* [33] in Latin. Its author was a German theologian and mathematician Johann Heilbronner (1706–1745/47). The book contained an index of names. It was a bulky book – 924 pages. Montucla called this Heilbronner's work "chaos" [10, p. 152-153]. The author paid much attention to philosophical issues and described the structure of mathematics. He presented the sequence of certain names and discoveries in mathematics in considerable detail, although not free from errors. He listed famous manuscripts and published books. This *Historia matheseos* compared favourably with the previous books thanks to the two special features. First, the author added modest information from the history of Arab and Chinese mathematics (names and discoveries) to the European history; second, all these different national histories were reduced to a single time scale. Heilbronner used the achievements in chronology of the last century and dated each event in mathematics in several ways: mentioned the eclipses which happened at that time or other celestial events, specifying their characteristics from astronomical tables (of Ptolemy and other astronomers), year anno mundi (ad annum Mundi), year from the founding of the City of Rome (ab urbe condita), year B.C. (ante Christum natum, ante Christi nativitatem), or year A.D. (ab Anno Christi). This presentation was not infallible. E.g. on page 353, Thang-Heng[38] (Zhang Heng, Chang Heng, 78–139), Chinese mathematician and astronomer, was dated to the year 164. Heilbronner carried Michael Psellos (11$^{th}$ century) back to the 9$^{th}$ century (p. 410), while Al-Farabi (872–950) and Ibn Musa (al-Khwarizmi, around 820) were dated to the 10$^{th}$ century. However, regardless of numerous disadvantages, it was in Heilbronner's book that the image of the history of world mathematics appeared for the first time ever, combining histories of different cultures. After Heilbronner's death, his library was purchased by Kästner who wrote a history of mathematics of his own. However, this is a topic for another article.

---

[35] E.g. epigram of Cl. Barksdale (1609–1687): "Savilian Bernard's a right learned man;/Josephus he will finish when he can".
[36] *Conics, a* fundamental treatise of Apollonius of Perga consisted of eight books. The Greek text of four of them has been preserved; three other books have survived translated into Arabic; the eighth book was renovated in the 18$^{th}$ century by E. Halley who published Apollonius' works (Oxford, 1710)
[37] Joseph Raphson, an English mathematician, Newton's disciple, died before 1715. We know very little about his life. He was the author of the most appropriate statement of Newton's approximation approach

The list of the books provided can be supplemented. A good review can be found in Popov's book [10], although one could blame him of some inaccuracies and gaps. However, the author was very scrupulous writing his book, he had read all the books he wrote about. A book devoted to historical development of historiography and mathematics in various countries was published in 2002 [34]. However, the period till 1750 was illustrated in it all too briefly.

Thus, in the first two millenia of its existence, the history of mathematics began developing scientific methodology: performing scientific analysis of works, sources (original, translated, renderings, and annotations); distinguishing between facts and interpretations; drafting catalogues and reference books; issues relating to individual and collective authorship (of national school), analysing the application and teaching of mathematical methods; chronology. In addition, textual analysis was emerging, the purpose of writing of mathematical works was not identified (research, teaching). Authors did not consider the role and reciprocal influence of ancient civilizations – as a rule, they began the history of mathematics from Greeks and considered it mostly in Latin culture. Arab manuscripts began to become the custom; Chinese manuscripts were hardly mentioned; and Indian manuscripts were almost unknown. The issues of national priorities were solved quite simply – each historian knew mathematical literature of his own country pretty well, giving preference to his compatriots (as, for example, Wallis or Wolff). This was the way the historical and mathematical memory of the nation consolidated and its mentality shaped. Absolute chronology was evolving up to the 17$^{th}$ century. Therefore, mathematical achievements in the social time of different civilisations only began to be compared. They did not distinguish between the development stages, periods of boom and bust, they did not give prominence to the trends in the evolution of mathematics, they did not emphasise its independence or the degree of its dependence on the needs of that time. Gradually, the presentation of the history of mathematics was undergoing stages – from describing discoveries and biographies by type of chronicles to genesis of ideas and understanding of mathematical progress.

*Conclusions.* As a science, the history of mathematics was evolving together with the evolvement of mathematics itself. Ancient works described the sequence of discoveries in mathematics based on the principle 'who discovered what', 'who taught whom', 'who followed whom'. Eudemus of Rhodes compared works that were topically related. Proclus began distinguishing research from educational works. With Dionysius de Laerte, biographies appeared in the history of mathematics. However, the history of mathematics was confined in a single culture.

---

[38] Heilbronner used to be very reliant on the letters about the history of Chinese astronomy of Antoine Gaubil, cartographer and missionary in China, which were published in European journals as of 1729.

The phenomenon of scientific translation as an art and its systematisation appeared in Arab culture of the 8th–12th centuries. Thanks to Arabic translations, the ancient heritage was preserved. This tradition was carried on in the Christian medieval period: Greek, Byzantine, and Arabic texts were translated into Latin; compendiums came into being – those were abridged renderings of classic authors which contained historical information. In Muslim culture of the 14th–17th centuries, great importance was attached to accounting and classifying manuscripts, describing thereof providing information from their authors' biographies; first catalogues began to be created. As of the 12th century, Englishmen began collecting manuscripts brought to England from the East; research libraries emerged, e.g. the Bodleian library (14th century). The advent of book printing (15th century) gave an impulse to dissemination of works of ancient classic authors; they were annotated and discussed, which included criticism as well. They began accounting each issue of a book as an individual scientific event. However, until the 16th century, the entire body of mathematical knowledge appeared to be static, lacking development (Ramus), although it was consistently establishing in time. Only the development of teaching methods was considered. It was the 17th century that a notion of progress in mathematics appeared for the first time in Voss' works. The history of creation gradually replaced the history of discoveries.

Chronologically, research works were arranged in an ordered fashion inside each culture: the chronology by Olympiad, since the creation of Rome; Muslim chronology – from Hegira; Chinese, by dynasty; Christian, Anno Domini, etc. Events, including mathematical events, which happened in various cultures were not related to one another; the story of each culture was stated independently. Wallis, for example, believed that Arabs knew algebra: "After Diophantus, not to say before him, the notion of powering was studied by Arab authors. However, they remained unaware of it in Europe for a rather long time." [18, p. 4-5]. This also resulted in "national shortsightedness", when discoveries of compatriots seemed closer and more important (Wallis, Wolff).

Thanks to the works of Scaliger and Petavius, the chronology reform of the 16th–17th centuries made it possible to bring historical events in compliance with astronomical phenomena and reduce them to a single scale which had a starting point and a direct reading-scale and a countdown-scale (A.D., B.C.). It should be noted that by this time, zero began to be apprehended as a reference point and a negative number, as a possibility of a countdown of steps, time, temperature (Wallis in the 17th century, Celsius in the 18th century). Authors began to include not only Christian mathematicians in the history of mathematics, Muslim mathematicians were mentioned as well (Baldi, Blankanus).

Methods used in other historical sciences began to be used in the history of mathematics. Those were methods of paleography, historical criticism, chronology, the doctrine of historical sources, evidence of their authenticity, ways to identify forgery, ancient written instruments and materials, styles. In the 18th century, name indices began to appear in books.

Authors began to try and present the history of mathematics not as a chronicle but as a history of ideas (Buteo in the 16th century, then Wallis in the 17th century).

All this made the history of mathematics ready for its next fledging period which began in 1758 when Montucla's History of Mathematics appeared.

# References

1. Прокл Диадох. Комментарий к первой книге «Начал» Евклида (перевод А.И. Щетникова). - Москва: Русский Фонд Содействия Образованию и Науке, 2013, 368 с. (Приложение № III к журналу «Аристей. Вестник классической филологии и античной истории»). Prokl Diadokh. Kommentarij k pervoj knige «Nachal» Evklida (perevod A.I.

Shchetnikova). - Moscow: Russkij Fond Sodejstviya Obrazovaniyu i Nauke, 2013, 368 s. (Prilozhenie № III k zhurnalu «Aristej. Vestnik klassicheskoj filologii i antichnoj istorii»)

2. Жмудь, Л. Я. История математики Евдема Родосского //Hyperboreus: St. Petersburg, 1997, № 3, с. 274—297. Zhmud' L. Ya. Istoriya matematiki Evdema Rodosskogo //Hyperboreus, 1997, 3, p. 274—297.

3. Les Prolégomènes d'Ibn Khaldoun. Troisième partie. Traduits en Français et commentés par William Mac Guckin, Baron de Slane, membre de l'Institut. Troisième partie des tomes XIX, XX et XXI des Notices et Extraits des Manuscrits de la Bibliothèque Nationale publiés par l'Institut de France (1863). Новое издание: Paris: Librairie orientaliste Paul Geuthner, 1938, 574 p., p. 94-127.

4. Katib Çelebi. Lexicon bibliographicum et encyclopaedicum, a Mustafa Ben Abdallah Katib Jelebi dicto et nomine Haji Khalfa celebrato compositum. Ad codicum Vindobonensium, Parisiensium, et Berolinensis fidem primum edidit Latine vertit et commentario indicibusque instruxit Gustavus Fluegel. In VII volumes. Leipzig: Flügel, Gustav Leberecht, 1835 – 1858.

5. Buteo, I. De quadratura circuli libri duo : vbi multorum quadraturae confutantur & ab omnium impugnatione defenditur Archimedes ; eiusdem Annotationum opuscula in errores Campani, Zamberti, Orontij, Peletarij, Io. Penae interpretum Euclidis. Lyon: apud Gulielmum Rovillium, 1559, 283 p.

6. Beckmann P. A History of Pi. New York: Golem Press, 1971/2015, 190 p.

7. Ramus, P. Proemium mathematicum in tres libros distributum. Paris: Wechel, 1567, 501 p.

8. Ramus, P. Scholarum mathematicarum libri unus et triginta. Basel: Eusebius Episcopius, 1569, 320 p.

9. Матвиевская Г. П. Рамус. 1515—1572. М.: Наука, 1981, 150 с. Matvievskaya G. P. Ramus. 1515—1572. Moskva: Nauka, 1981, 150 p.

10. Попов, Г.Н. История математики. Москва Типо-Лит. Московского Картоиздательского Отдела корп. Воен. Топогр., 1920, 236 с. Popov G.N. Istoriya matematiki. Moskva, 1920, 236 p.

11. Baldi, B. Cronica de' matematici: overo Epitome dell' istoria delle vite loro. Urbino: Angelo Antonio Monticelli, 1707, 156 p.

12. Blancanus, J. De mathematicarum Natura dissertatio. Una cum Clarorum mathematicorum chronologia. Bologna: apud Bartholomaeum Cochium, 1615, 53 p.

13. Vossius, G.J. De universae matheseos natura et constitutione liber; cui subjungitur chronologia mathematicorum. Amsterdam: ex typogr. J. Blaeu, 1650, 473+32 p.

14. Vossius, G.J. De quatuor artibus popularibus, de philologia, et scientiis mathematicis. Cui operi subjungitur Chronologia Mathematicorum. Libri tres. Amsterdam: ex typographeio Ioannis Blaeu, 1660, 467+ 35 p.

15. Chales (Dechales), C.-F. Cursus seu Mundus Mathematicus. In 3 volumes. Lyon: ex officina Annisoniana, 1674. V.1 – 763 p., V.2 -731 p., V.3 – 863 p.

16. Smith, D.E. History of Mathematics, Vol. 1. New York: Dover Publications Inc., 1951, 613 p.

17. Mabillon, J. De re diplomatica, libri VI. Paris: Louis Billaine, 1681, 634 p.

18. Wallis, J. Treatise of algebra both Historical and Practical. Shewing, the Original, Progress and Advancement thereof, from time to time; and by what Steps it hath attained to the Heighth at which now it is. With some additional Treatises. London: Richard Davis. M.DC.LXXXV (1685), 374 p.

19. Speidell, J. New logarithmes: the first inuention whereof, was, by
the honourable Lo. Iohn Nepair, Baron of Marchiston, and printed at
Edinburg in Scotland, anno 1614, in whose vse was and is required the
knowledge of algebraicall addition and substraction, according to + and −. London, 1619. 32 p.

20. Speidell, J. New Logarithmes. London, 1622. [Reprinted in: Francis Maseres. Scriptores logarithmici, volume 6. London: R. Wilks, 1807, pp. 728-759].

21. Hobson, E. W. John Napier and the Invention of Logarithms, 1614: A Lecture by E.W. Hobson. Cambridge University Press, 1914/2012, 48 p.

22. Oughtred, G. Clavis Mathematicae. Oxoniae: typis Lichfieldianis, 1631; Key of the Mathematics (in English). London: Salusburn, 1694, 208 p.

23. Harriot, T. Artis Analyticae Praxis ad Aequationes Algebraicas Resolvendas. London: apud Robertum Barker, 1631. 180 p. [Thomas Harriot's Artis analyticae praxis: an English translation with commentary M. Seltman, R. Goulding, editors and translators. New York: Springer, 2007, 299 p.].

24. Эйлер, Л. Дифференциальное исчисление. Москва-Ленинград: Гос. изд-во технико-технической литературы, 1949, 580 с. Euler, L. Differencial'noe ischislenie. Moskva-Leningrad: Gos. izd-vo tekhniko-tekhnicheskoj literatury, 1949, 580 p.

25. Декарт, Р. Геометрия. Перевод, примечания и статья А.П. Юшкевича. Москва-Ленинград:Научтехиздат, 1938, 297 с. Descartes R. Geometriya. Perevod, primechaniya i stat'ya A.P. Yushkevicha. Moscow-Leningrad:Nauchtekhizdat, 1938, 297 p.

26. Вилейтнер, Г. Истории математики от Декарта до середины XIX столетия. М.: ГИФМЛ, 1960, 468 с. Wieleitner, H. Istorii matematiki ot Dekarta do serediny XIX stoletiya. M.: Gosudarstvennoe izdatel'stvo fiziko-matematicheskoj literatury, 1960, 468 p.

27. Bernard, E. Catalogi librorum manuscriptorum Angliae et Hiberniae in unum collecti cum indice alphabetico. 2 volumes. Oxford: e Theatro Sheldoniano, 1697. 1011+76 p.

28. Robert Huntington (bp. of Raphoe), Edward Bernard, Thomas Smith. Epistolae: Et Veterum Mathematicorum, Graecorum, Latinorum, & Arabum, Synopsis. London: typis G. Bowyer, impensis A. & J. Churchill, 1704. 2 v. part 9, p.1-78.

29. Bernard, E. Et Veterum Mathematicorum, Graecorum, Latinorum, & Arabum, Synopsis// Robert Huntington (bp. of Raphoe), Edward Bernard, Thomas Smith. Epistolae: Et Veterum Mathematicorum, Graecorum, Latinorum, & Arabum, Synopsis. London: Typis G. Bowyer, impensis A. & J. Churchill, 1704, part 11, p. 1-44.

30. Raphson, J. Historia fluxionum : sive tractatus originem & progressum peregregiæ istius methodi brevissimo compendio (et quasi synopticè) exhibens. Per Josephum Raphsonum (History of Fluxions, Showing in a Compendious Manner the First Rise of, and various Improvements made in that Incomparable method defending Newton's role). London: printed by William Pearson, 1715, 123 p.

31. Wolff, Chr. Vollständiges mathematisches Lexicon, darinnen alle Kunst-Wörter und Sachen, welche in der erwegenden und ausübenden Mathesi vorzukommen pflegen, deutlich erkläret; überall aber zur Historie der mathematischen Wissenschafften dienliche Nachrichten eingestreuet, und die besten und auserlesensten Schrifften, welche jede Materie gründlich abgehandelt, angeführet : ferner auch die Mund- und Redens-Arten derer Marckscheider auch hieher gehöriger Künstler und Handwercker beschrieben; und endlich alles zum Nutzen so wohl gelehrter als ungelehrter liebhaber der vortrefflichen mathematick eingerichtet worden. Nebst XXXVI. kupffer-tabellen, Leipzig: Gleditsch, 1716, 1734, 1747, 788 c.

32. Wolfius, Chr. Elementa matheseos Universae. Tomus V. Halle: Renger, 1741. 340 p.

33. Heilbronner, J.C. Historia matheseos universae a mundo condito ad saeculum post Christ. nat. XVI. Accedit recensio elementorum, compendiorum et operum mathematicorum atque historia arithmetices ad nostra tempora. Leipzig: Gleditsch, 1742. 924 c.

34. Writing the History of Mathematics: Its Historical Development. Ed. J. W. Dauben, Christoph J. Scriba. Birkhäuser, Basel, 2002, 689 c.

35. Синкевич Г.И. Ранний этап развития historia matheseos. Историография истории математики // История науки и техники, 2017 г. №1. С. 3-17.